\documentclass[11pt]{article}

\usepackage{color}
\usepackage{latexsym}
\usepackage{amssymb}
\usepackage{graphicx}
\usepackage{amsmath, amsfonts,amssymb,theorem,euscript,array,enumerate,amsfonts,mathrsfs}

\newtheorem{Remark}{Remark}[part]

\def\esssup_#1{\underset{#1}{\mathrm{ess\,sup\, }}}
\def\essinf_#1{\underset{#1}{\mathrm{ess\,inf\, }}}
\def\argmax_#1{\underset{#1}{\mathrm{arg\,max\, }}}

\def \Frac{\displaystyle\frac}

\def \trans{^{\scriptscriptstyle{\intercal }}}

\def \N{\mathbb{N}}
\def \R{\mathbb{R}}

\def \E{\mathbb{E}}
\def \F{\mathbb{F}}
\def \G{\mathbb{G}}
\def \H{\mathbb{H}}
\def \L{\mathbb{L}}
\def \P{\mathbb{P}}

\def \S{\mathbb{S}}

\def \Ac{{\cal A}}
\def \Bc{{\cal B}}

\def \Dc{{\cal D}}
\def \Ec{{\cal E}}
\def \Fc{{\cal F}}
\def \Gc{{\cal G}}

\def \Lc{{\cal L}}
\def \Pc{{\cal P}}
\def \Mc{{\cal M}}

\def \Uc{{\cal U}}

\def \Yc{{\cal Y}}

\def \Zc{{\cal Z}}

\def \ep{\hbox{ }\hfill$\Box$}

\def\Dt#1{\Frac{\partial #1}{\partial t}}

\def\reff#1{{\rm(\ref{#1})}}

\def\beqs{\begin{eqnarray*}}
\def\enqs{\end{eqnarray*}}
\def\beq{\begin{eqnarray}}
\def\enq{\end{eqnarray}}

\begin{document}

\title{Feynman-Kac representation of fully nonlinear PDEs \\ and applications}

\author{Huy\^en PHAM
             \\\small  Laboratoire de Probabilit\'es et  Mod\`eles Al\'eatoires
             \\\small  CNRS, UMR 7599
             \\\small  Universit\'e Paris Diderot, Paris, France 
             \\\small  pham at math.univ-paris-diderot.fr
             \\\small  and John Von Neumann Institute, Ho-Chi-Minh City 
             }

%\date{}

\maketitle

\begin{abstract}
The classical Feynman-Kac formula  states the connection between linear parabolic partial differential equations (PDEs), like  the heat equation, and expectation of stochastic processes driven by Brownian motion. 
It gives then a method for solving linear PDEs by Monte Carlo simulations of random processes. The  extension to (fully)nonlinear PDEs  led in  the recent years to important developments in stochastic analysis and the emergence of the theory of backward stochastic differential equations (BSDEs), which can be viewed as nonlinear Feynman-Kac  formulas.  We review in this paper the main ideas and results in this area, and present implications of these probabilistic representations for the numerical resolution of nonlinear PDEs, together with some applications to  stochastic control problems and model uncertainty  in finance. 
 \end{abstract}

\vspace{5mm}

\noindent {\bf MSC Classification (2000):}  60H30, 65C99, 93E20

\vspace{5mm}

\noindent {\bf Keywords:}   Nonlinear PDE, Hamilton-Jacobi-Bellman equation, backward stochastic diffe\-rential equation, randomization of controls, discrete-time approximation.

\newpage

\section{Introduction}

 \setcounter{equation}{0} \setcounter{Assumption}{0}
\setcounter{Theorem}{0} \setcounter{Proposition}{0}
\setcounter{Corollary}{0} \setcounter{Lemma}{0}
\setcounter{Definition}{0} \setcounter{Remark}{0}

Let us consider the parabolic heat equation:

\begin{equation} \label{heat}
\left\{
\begin{array}{rccl}
\Dt{v} + \frac{1}{2} \Delta_x  v &=& 0, &   \mbox{ on } [0,T)\times\R^d, \\
%& & \\
v(T,.) &=& h, &  \mbox{ on } \;  \R^d.
\end{array}
\right.
\end{equation}
It is well-known that the solution to \reff{heat} is given by:
\beqs
v(t,x) &=&  \int h(y) K(T-t,x,y) dy, 
\enqs
where $K(t,x,y)$ $=$ $\frac{1}{(4\pi t)^{d\over 2}} e^{-|x-y|^2/4t}$ is  the heat kernel on $\R^d$.  By introducing the $d$-dimensional Brownian $W$ on a probability space $(\Omega,\Fc,\P)$, and from the Gaussian distribution 
of  $W_t$, we observe that the solution $v$ can be represented also as:
\beq \label{feynheat}
v(t,x) &=& \E\big[ h(x  + W_{T-t}) \big], \;\;\; (t,x) \in [0,T]\times\R^d. 
\enq
The probabilistic representation \reff{feynheat} gives a Monte-Carlo method for computing an appro\-ximation of $v$ by the empirical mean:
\beqs
v(t,x) & \simeq & \bar v^N(t,x)  \; := \;  \frac{1}{N} \sum_{i=1}^N h(x + W_{T-t}^i), 
\enqs
where $(W^i)_{1\leq i\leq N}$ is an $N$-sample drawn from an (exact) simulation of $W$.  The convergence of $\bar v^N$  to $v$ is ensured by the law of large numbers, when $N$ goes to infinity, while the rate 
of convergence, obtained from the central limit theorem, is equal to $1/\sqrt{N}$, and independent of the dimension $d$ of the heat equation.  More generally,  let us consider the linear parabolic partial differential equation (PDE):
\begin{equation} \label{linearPDE}
\left\{
\begin{array}{rccl}
\Dt{v} + \Lc  v + f  &=& 0, &   \mbox{ on } [0,T)\times\R^d, \\
%& & \\
v(T,.) &=& h, &  \mbox{ on } \;  \R^d,
\end{array}
\right.
\end{equation}
 where $\Lc$ is  the second order Dynkin operator:
 \beq \label{dynkin}
 \Lc v &=& b(x). D_x v + \frac{1}{2}{\rm tr}(\sigma\sigma\trans(x) D_x^2 v). 
 \enq
Under suitable conditions on the functions $b$, $\sigma$, $f$ and $h$ defined on $\R^d$, there exists a unique solution $v$ to \reff{linearPDE}, which may be represented by the Feynman-Kac formula: 
\beq \label{feynlin}
v(t,x) &=& \E \Big[ \int_t^T f(X_s^{t,x}) ds + h(X_T^{t,x}) \Big], \;\;\; (t,x) \in [0,T]\times\R^d, 
\enq
where $X^{t,x}$ is the solution to the (forward) diffusion process, 
\beqs
dX_s &=& b(X_s) ds + \sigma(X_s) dW_s,  \;\;\;  s \geq t, 
\enqs
starting from $x$ at time $t$.  
Notice that the Feynman-Kac formula \reff{feynlin} can be easily derived from It\^o's formula when $v$ is smooth. Indeed, in this case, by defining the pair of processes 
$(Y,Z)$:
\beqs
Y_t \; := \; v(t,X_t), & & Z_t \; := \; \sigma\trans(X_t) D_x v(t,X_t), \;\;\; 0 \leq t \leq T,
\enqs
and applying It\^o's formula to $v(s,X_s)$ between $t$ and $T$,  with $v$ satisfying the PDE \reff{linearPDE}, we get: 
\beq \label{BSE}
Y_t &=& h(X_T) + \int_t^T f(X_s) ds - \int_t^T Z_s dW_s, \;\;\; 0 \leq t \leq T. 
\enq
This equation can be viewed as a backward stochastic equation in the pair of adapted processes $(Y,Z)$ w.r.t. the filtration $\F^W$ generated by the Brownian motion $W$, 
determined from a terminal condition $h(X_T)$, and originally appeared in \cite{bis73}. 
By taking conditional expectation in \reff{BSE}, we retrieve the Feynman-Kac formula 
\reff{feynlin}.  This probabilistic representation leads to a numerical method for solving  the linear PDE, relying on Monte-Carlo simulations of the forward diffusion process $X$, whose convergence rate does not depend on the dimension of the problem, hence not suffering in principle of the curse of dimen\-sionality encountered in deterministic numerical methods. 
On the other hand, it is also useful for computing explicitly the solution $v$ in some particular models for $X$, e.g.  geometric Brownian motion  
in the Black Scholes model for option pricing in finance.

In this paper, we address the  problem  of nonlinear PDEs, and shall review the recent developments about their probabilistic representation, i.e. nonlinear 
Feynman-Kac formulae. We  shall first consider in Section 2
 the case of semi-linear PDEs, i.e. when nonlinearity appears only on the first order derivative, and show how it is related to the 
theory of backward stochastic differential equations (BSDEs)  introduced in \cite{parpen90}, and leads to probabilistic scheme for solving semi-linear PDEs. 
 We next consider in Section 3  the challenging problem of fully nonlinear PDEs, i.e. when nonlinearity enters also on the second order derivative. Such framework arises  in many applications, for example in stochastic control in finance (portfolio optimization, risk management, model uncertainty). We  shall  present  the randomization approach for dealing with such nonlinear context, and show how fully nonlinear PDES  are represented in terms of randomized BSDEs with nonpositive jumps. This provides an original probabilistic scheme for solving fully nonlinear PDEs.

\section{Backward SDEs and semi-linear PDEs }

 \setcounter{equation}{0} \setcounter{Assumption}{0}
\setcounter{Theorem}{0} \setcounter{Proposition}{0}
\setcounter{Corollary}{0} \setcounter{Lemma}{0}
\setcounter{Definition}{0} \setcounter{Remark}{0}

\subsection{A short overview of BSDEs}

Let us introduce some standard notations in the theory of backward stochastic differential equations (BSDEs). On a complete probability space 
$(\Omega,\Fc,\P)$ on which is defined a $d$-dimensional Brownian motion $W$ over a finite time interval $[0,T]$,  
and its natural filtration $\F$ $=$ $\F^W$, we denote by: 
\begin{itemize}
\item $\Pc_\F$: $\sigma$-algebra of $\F$-predictable subsets of  $[0,T]\times\Omega$
\item $\S_\F^2$:  set of  real-valued c\`ad-l\`ag $\F$-adapted processes $Y$  such that 
\beqs
\E\big[\sup_{0\leq t\leq T}  |Y_t|^2 \big] &<& \infty,
\enqs
\item $\L_\F^2(W)$: set of  $\R^d$-valued $\Pc_\F$-measurable processes $Z$ such that 
\beqs
\E\left[\int_0^T |Z_t|^2dt\right] &<& \infty.
\enqs
\end{itemize}

We are given as data:
\begin{itemize}
\item a terminal condition $\xi$, which is an $\Fc_T$-measurable real-valued random variable 
\item a generator $f$ $=$ $(f_t(y,z))_{0\leq t\leq T}$, which is an $\Pc_\F\otimes\Bc(\R\times\R^d)$-measurable real-valued map, where $\Bc(\R\times\R^d)$ denotes the Borel $\sigma$-field of $\R\times\R^d$.
\end{itemize}
A (one dimensional) BSDE in differential form is written as
\beq \label{defBSDE}
dY_t &=& - f_t(Y_t,Z_t) dt + Z_t dW_t, \;\;\; 0 \leq t \leq T, \;\; Y_T = \xi,
\enq
and a solution to \reff{defBSDE} is a pair $(Y,Z)$ $\in$ $\S_\F^2\times\L_\F^2(W)$ satisfying
\beq \label{BSDE}
Y_t &=& \xi + \int_t^T  f_s(Y_s,Z_s) ds - \int_t^T Z_s  dW_s, \;\;\; 0\leq t \leq T.
\enq

Existence and uniqueness of a solution to the BSDE \reff{defBSDE} is proved in the seminal paper \cite{parpen90} under the following Lipschitz and square integrability assumptions: 

\vspace{2mm}

{\bf (H1)} 
\begin{itemize}
\item[(i)]  $f$ is uniformly Lipschitz in $(y,z)$, i.e. there exists a positive constant $C_f$ s.t. for all $(y,z,y',z')$: 
\beqs
|f_t(y,z) - f_t(y',z')| & \leq & C_f \Big( |y-y'| + |z-z'|\Big), \;\; dt\otimes d\P \; a.e. 
\enqs
\item[(ii)] $\xi$ and $\{f_t(0,0),t\in [0,T]\}$ are square integrable: 
\beqs
\E \Big[ |\xi|^2 + \int_0^T |f_t(0,0)|^2 dt \Big] & < & \infty. 
\enqs
\end{itemize}
Notice that  when the generator $f$ does not depend on $y$ and $z$, the solution to the BSDE \reff{defBSDE}, which is then simply a backward stochastic equation as in \cite{bis73},  is directly obtained from the martingale representation theorem applied to the 
Brownian martingale $M_t$ $:=$ $\E \Big[ \xi + \int_0^T f_s ds \big | \Fc_t]$, $0\leq t\leq T$, which gives the existence of an integrand $Z$ $\in$ $\L_\F^{2}(W)$ s.t. 
\beqs
M_t & = & M_0 + \int_0^t Z_s dW_s , \;\;\; 0 \leq t \leq T.
\enqs
Indeed,  by defining 
\beqs
Y_t := M_t-\int_0^t f_s ds &=& \E \Big[ \xi + \int_t^T f_s ds \big | \Fc_t] , \;\;\; 0\leq t\leq T,
\enqs
 we see that $(Y,Z)$ satisfies \reff{BSDE}. In the general case where $f$ depends on $(y,z)$, the existence and uniqueness is proved by a fixed point argument under the Lipschitz assumption in {\bf (H1)}(i).

\vspace{2mm}

Let us now consider the particular case of interest when the generator is linear, i.e. in the form:
\beqs
f_t(y,z) &=& \delta_t y + \alpha_t. z + \gamma_t, 
\enqs
for some bounded $\F$-adapted processes $(\delta_t)$ valued in $\R$, $(\alpha_t)$ valued in $\R^d$, and 
$(\gamma_t)$ $\in$ $\H_\F^2$ the set of rela-valued $\F$-adapted processes s.t. $\E[\int_0^T |\gamma_t|^2 dt]$ $<$ $\infty$. 
By discounting and Girsanov's change of measure, the solution (in $Y$) to the linear BSDE \reff{BSDE} is given by the linear expectation: 
\beqs
Y_t &=& \E^{\P^\alpha}\big[ e^{\int_t^T \delta_s ds} \xi + \int_t^T  e^{\int_t^u \delta_u du} \gamma_s ds \big| \Fc_t \big],  
\enqs
where $\P^\alpha$ is the probability measure equivalent to $\P$ under which 
\beq \label{defWalpha}
W^\alpha  &:=& W - \int \alpha dt, \;\;\; \mbox{ is a } \; \P^\alpha-\mbox{Brownian motion}. 
\enq
Such context of linear BSDE arises typically in option pricing in finance, where $\P^\alpha$  is the martingale measure, $Y$ is the fair price for the option payoff $\xi$,  and $Z$ the hedging portfolio.  
A trivial remark in that, in this linear case, if  $\xi$ $\geq$ $0$ and $\gamma$ $\geq$ $0$, then the  solution $Y$ to the linear BSDE is also nonnegative. This is the key observation for showing comparison theorem for BSDEs:  given two pairs $(\xi,f)$ and $(\xi',f')$ of terminal data/generators satisfying {\bf (H1)}, and $(Y,Z)$, $(Y',Z')$ be the solutions to their BSDEs.  Suppose that:
\beqs
\xi \; \leq \; \xi'  \; a.s.  & \mbox{ and } & f_t(Y_t,Z_t) \; \leq \; f'_t(Y_t,Z_t), \;\; dt\otimes d\P \; a.e. 
\enqs
Then, 
\beqs
Y_t & \leq & Y_t', \;\;\; 0 \leq t \leq T. 
\enqs

Let us next consider another case of interest where the generator $f_t(y,z)$ is convex in $z$, and is written in the form:
\beq \label{convexf}
f_t(y,z) &=& \delta_t y + \sup_{a \in A}\big[a.z + \gamma_t(a)\big], 
\enq
for some bounded adapted processes $(\delta_t)$,  where $A$ is a compact subset of $\R^d$, $\gamma_t(a)$ is a $\Pc_\F\otimes\Bc(A)$-measurable map s.t. 
$\esssup_{a\in A}|\gamma_t(a)|$ $\in$ $\H_\F^2$. By using comparison theorem for BSDEs and result for the linear case, one shows that the solution (in $Y$) to the BSDE is expressed as: 
\beqs
Y_t &=& \esssup_{\alpha \in \Ac} \E^{\P^\alpha} \Big[ e^{\int_t^T \delta_s ds} \xi + \int_t^T e^{\int_t^u \delta_u du} \gamma_s(\alpha_s) ds  \big| \Fc_t \Big], \;\;\; 0 \leq t \leq T, 
\enqs
where $\Ac$ is the set of adapted processes $\alpha$ valued in $A$ and $\P^\alpha$ is the probability measure equivalent to $\P$ under which the drifted process 
$W^\alpha$ defined in \reff{defWalpha} is a  $\P^\alpha$-Brownian motion.  Hence, controlled drift problems and risk measures with uncertain drifts 
are  related to BSDE by choosing a generator in the form \reff{convexf}.  We refer to \cite{elkmaz97}, \cite{elkpenque97} or \cite{pen03} for a more detailed review and applications of BSDEs.

\subsection{Markov BSDEs and PDEs}

We put ourselves in a Markov setting in the sense that we suppose that  the terminal data and generator of the BSDE \reff{BSDE} are in the form: 
\beqs
\xi \; = \; h(X_T), & & f_t(\omega,y,z) \; = \; f(X_t(\omega),y,z)
\enqs
where $h(x)$ is some measurable function on $\R^d$,  $f(x,y,z)$ is some measurable function on $\R^d\times\R\times\R^d$ (we kept the same notation $f$ by misuse), and 
$X$ is  a forward  diffusion process of dynamics: 
\beq \label{defX}
dX_s &=& b(X_s) ds + \sigma(X_s) dW_s \;\;\; \mbox{ in } \; \R^d.  
\enq
Under standard Lipschitz assumptions on the coefficients $b$ $:$ $\R^d$ $\mapsto$ $\R^d$, and $\sigma$ $:$ $\R^d$ $\mapsto$ $\R^{d\times d}$,  there exists a unique strong solution to 
\reff{defX}  given some initial condition, and we have the standard estimate:
\beqs
\E \big[ \sup_{0\leq t\leq T} | X_t|^2 \big]  & \leq & C(1 + |X_0|^2). 
\enqs 
A forward BSDE is then written as:
\beq \label{BSDEmarkov}
Y_t &=& h(X_T)  + \int_t^T  f(X_s,Y_s,Z_s) ds - \int_t^T Z_s  dW_s, \;\;\; 0\leq t \leq T,
\enq
and under Lipschitz condition on $f$, and linear growth condition on $h$, there exists a unique solution $(Y,Z)$ to the Markov BSDE \reff{BSDEmarkov}.  Denoting by  $(Y^{t,x},Z^{t,x})_{t\leq s\leq T}$ the solution  to the BSDE \reff{BSDEmarkov} when $X$ $=$ $(X^{t,s})_{t\leq s\leq T}$ is the solution to \reff{defX} starting from $x$ at time $t$, we notice that 
\beq \label{defvBSDE}
v(t,x) &:=& Y_t^{t,x}, \;\;\; (t,x) \in [0,T]\times\R^d, 
\enq
is a deterministic function on  $[0,T]\times\R^d$, and by the Markov property of the diffusion process, we have:
\beqs
Y_t &=&  v(t,X_t), \;\;\; 0 \leq t \leq T. 
\enqs
 
Let us now derive formally the PDE satisfied by the function $v$.   By definition of the Markov BSDE \reff{BSDEmarkov}, we have:
 \beqs
 Y_s - Y_t \; = \; v(s,X_s) - v(t,X_t) &=& -  \int_t^s f(X_u,Y_u,Z_u) ds  +  \int_t^s Z_u dW_u,
 \enqs
for all $0\leq t\leq s\leq T$.  Assuming that $v$ is smooth, it follows from It\^o's formula:
 \beqs
 & & \int_t^s (\Dt{v} + \Lc v)(u,X_u)  du + \int^s \sigma\trans(X_u)D_x v(u,X_u) dW_u \\
 &=& -  \int_t^s f(X_u,Y_u,Z_u)  du  +  \int_t^s Z_u dW_u,
 \enqs
where $\Lc$ is the Dynkin operator associated to the diffusion $X$, and given in \reff{dynkin}.  Identifying  the   finite variation terms in ``$dt$" and the Brownian martingale terms in  ``$dW$", 
we see that  
\beqs
Z_t &=& \sigma\trans(X_t)D_x v(t,X_t)
\enqs
and $v$ should satisfy the semi-linear parabolic PDE:
\begin{equation} \label{semilinearPDE}
\left\{
\begin{array}{rccl}
\Dt{v} + \Lc  v + f(x,v,\sigma\trans D_x v)  &=& 0, &   \mbox{ on } [0,T)\times\R^d, \\
%& & \\
v(T,.) &=& h, &  \mbox{ on } \;  \R^d. 
\end{array}
\right.
\end{equation}
The main issue in this derivation  comes from the fact that in general, the function $v$ is not  smooth, and this is overcome with the notion of viscosity solution: it is proved in \cite{parpen92} that the function $v$ 
in \reff{defvBSDE} is the unique viscosity solution to \reff{semilinearPDE}.  Therefore,  the BSDE \reff{BSDEmarkov}  provides a probabilistic representation to the solution of the semi-linear PDE \reff{semilinearPDE}.  This  extends the  Feynman-Kac formula \reff{feynlin} to the case where $f$ $=$ $f(x,y,z)$ depends on $y,z$, and we shall see in the next paragraph how it can be used to design a probabilistic numerical scheme for computing the solution to the semi-linear PDE  \reff{semilinearPDE}.

\subsection{Numerical issues}

The first step in the numerical scheme for the resolution of the BSDE \reff{BSDEmarkov} is the discrete-time approximation. It is constructed as follows. 

\vspace{1mm}

\noindent $\bullet$ {\it Euler scheme for the forward process}.  We are given a time grid $\pi$ $:=$ $\{t_0=0 < t_1 < \ldots < t_n=T\}$ of $[0,T]$, with modulus $|\pi|$ 
$:=$ $\max_{i=1,\ldots,n} \Delta t_i$, $\Delta t_i$ $:=$ $t_{i+1}-t_i$,  and approximate the forward diffusion process $X$ by its Euler scheme $X^\pi$ defined as
\beqs
 X^{\pi}_{t_{i+1}} &:=&  X^{\pi}_{t_i} + b(X^{\pi}_{t_i}) \Delta t_i + \sigma(X^{\pi}_{t_i}) \Delta W_{t_i}, \;\;\; i < n, \;\;\; X_{0}^\pi \; = \; X_0, 
 \enqs
where $\Delta W_{t_i}$ $:=$ $W_{t_{i+1}}-W_{t_i}$.   

\vspace{1mm}

\noindent $\bullet$ {\it Euler scheme for the backward process}. We first approximate the terminal condition $Y_T$ $=$ $h(X_T)$ by simply replacing $X$ by its Euler scheme: 
$Y_T$ $\simeq$ $h(X_T^\pi)$. Then, from the  formal Euler backward discretization:
\beqs
Y_{t_i} &=& Y_{t_{i+1}} + \int_{t_i}^{t_{i+1}} f(X_s,Y_s,Z_s) ds - \int_{t_i}^{t_{i+1}} Z_s dW_s  \\
& \simeq & Y_{t_{i+1}}  +  f(X_{t_i}^\pi,Y_{t_i},Z_{t_i}) \Delta t_i -  Z_{t_i} \Delta W_{t_i}. 
\enqs
we define the discrete-time approximation of the BSDE as follows:
\begin{itemize}
\item[(1)] taking expectation conditionally on $\Fc_{t_i}$ on both sides yields
\beqs
Y_{t_i} & \simeq & \E \big[ Y_{t_{i+1}} | \Fc_{t_i} \big] + f(X_{t_i}^\pi,Y_{t_i},Z_{t_i}) \Delta t_i 
\enqs
\item[(2)] Multiplying by $\Delta W_{t_i}$ and then taking conditional expectation gives
\beqs
0 & \simeq & \E \big[ Y_{t_{i+1}} \Delta W_{t_i} | \Fc_{t_i} \big] -  Z_{t_i} \Delta t_i.
\enqs
\end{itemize}
This formal approximation argument leads to a backward Euler scheme $(Y^\pi,Z^\pi)$ of the form: 
\begin{equation} \label{Eulerbackward}
\left\{
\begin{array}{ccl}
Z_{t_i}^\pi &=& \E \Big[ Y_{t_{i+1}}^\pi \frac{\Delta W_{t_i}}{\Delta t_i}   \big | \Fc_{t_i}  \Big],  \\
Y_{t_i}^\pi &=&  \E \big[ Y_{t_{i+1}}^\pi | \Fc_{t_i} \big] + f(X_{t_i}^\pi,Y_{t_i}^\pi,Z_{t_i}^\pi) \Delta t_i, \;\;\; i < n,  
\end{array}
\right.
\end{equation}
with terminal condition $Y_{t_n}^\pi$ $=$ $h(X_{t_n}^\pi)$.

\begin{Remark}
{\rm  The above scheme is implicit as $Y_{t_i}^\pi$ appears in both sides of the equation. Since $f$ is assumed to be Lipschitz and since it is multiplied by 
$\Delta t_i$, intended to be small, the equation can be solved numerically very quickly by standard fixed point methods. Alternatively, we could also consider an explicit scheme by replacing the second equation 
in \reff{Eulerbackward} by
\beqs
Y_{t_i}^\pi &=&  \E \big[ Y_{t_{i+1}}^\pi +  f(X_{t_i}^\pi,Y_{t_{i+1}}^\pi,Z_{t_i}^\pi) \Delta t_i  | \Fc_{t_i} \big].
\enqs
This will not change the convergence rate. 
\ep
}
\end{Remark}

\vspace{2mm}

The discrete-time approximation is measured by the squared error: 
\beqs
\Ec(\pi)^2 & :=&  \max_{i \leq n} \E \big[ |Y_{t_i}-Y_{t_i}^\pi|^2 \big]  +  \sum_{i=0}^{n-1} \E \Big[ \int_{t_i}^{t_{i+1}} |Z_t - Z_{t_i}^\pi|^2 dt \Big]. 
\enqs
By using It\^o's formula, Gronwall's lemma and Young inequality, it is proved in \cite{boutou04} and \cite{zha04} that under the Lipschitz assumption on the driver $f$, there exists a 
constant $C$ independent of $\pi$ such that: 
\beq
\Ec(\pi)^2  & \leq & C \Big(  \E\big|h(X_T) - h(X_T^\pi)|^2 + \max_{i\leq n} \E\big| X_{t_i} - X_{t_i}^\pi\big|^2  \nonumber \\
& & \;\;\;\;\;\;  + \;    \sum_{i=0}^{n-1} \E \Big[ \int_{t_i}^{t_{i+1}} |Z_t - \bar Z_{t_i}|^2 dt \Big] \Big),  \label{decerror}
\enq
where $\bar Z_{t_i}$ $=$ $\frac{1}{\Delta t_i} \E\big[ \int_{t_i}^{t_{i+1}} Z_t dt | \Fc_{t_i} \big]$.  In other words, we have three different error contributions:
\begin{itemize}
\item[1.] Strong approximation of the terminal condition, which depends on the terminal data and the forward Euler scheme 
\item[2.] Strong approximation of the forward SDE, which depends on the forward Euler scheme, but not on the BSDE problem
\item[3.] $L^2$-regularity of $Z$, which is intrinsic to the BSDE problem. 
\end{itemize}
It is well-known, see e.g. \cite{klopla92}, that the strong approximation of the forward Euler scheme:  $ \max_{i\leq n} \E\big| X_{t_i} - X_{t_i}^\pi\big|^2$ provides an error of order $|\pi|$. 
Consequently, when the terminal data $h$ is Lipschitz, this also  gives an error for $\E\big|h(X_T) - h(X_T^\pi)|^2$ of order $|\pi|$.  Finally,  it is proved in \cite{zha04} that under Lipschitz condition on $f$, we have the $L^2$-regularity of $Z$:  $\sum_{i=0}^{n-1} \E \big[ \int_{t_i}^{t_{i+1}} |Z_t - Z_{t_i}^\pi|^2 dt \big]$ $=$ $O(|\pi|)$. Therefore, under Lipschiz assumption on $f$ and  $h$, the rate of convergence of  the discrete-time approximation error  $\Ec(\pi)$ is of order  $|\pi|^{1\over 2}$:  
\beqs
\Ec(\pi) & \leq & C |\pi|^{1\over 2}. 
\enqs
This convergence rate is clearly optimal and similar to the one obtained for forward SDEs.

\vspace{2mm}

The practical implementation of the numerical scheme \reff{Eulerbackward} requires the computation of  conditional expectations. The key observation in our Markovian context is that all these conditional expectations are regressions, i.e. 
\beqs
\E \big[ Y_{t_{i+1}}^\pi   | \Fc_{t_i}  \big] \; = \; \E \big[  Y_{t_{i+1}}^\pi  | X_{t_i}^\pi \big], & & 
\E \big[ Y_{t_{i+1}}^\pi \Delta W_{t_i}    | \Fc_{t_i}  \big] \; = \;  \E \big[  Y_{t_{i+1}}^\pi \Delta W_{t_i}   | X_{t_i}^\pi \big], 
\enqs
which can be approximated by methods from statistics: 
\begin{itemize}
\item {\it Quantization.} Each $X_{t_i}^\pi$ is replaced by a quantized version, i.e. a projection on a finite grid, which is computed in some optimal way based on stochastic algorithm and Monte-Carlo simulations of $X^\pi$ (Kohonen).  The conditional expectation is then reduced to a discrete sum with weights also computed off line as the grid points.   We refer to \cite{balpag03}, 
\cite{pagphapri04} and the references therein. 
\item {\it Integration by parts.} The conditional expectation is approximated via an integration by parts  formula and Malliavin calculus, see \cite{boueketou04}. 
\item {\it  Least-square regression.}  The conditional expectation is approximated by non-parame\-tric regression methods, and the most popular one, known as Longstaff-Schwartz method  
\cite{lonsch01}, consists in  
the projection  on a set of basis functions, with optimal coefficients  computed from empirical least-square based on Monte-Carlo simulations of $X_{t_i}^\pi$.  We refer to 
 \cite{gobetal06}  for more details and analysis of convergence rate of this approach. 
\end{itemize}
The advantage of these probabilistic methods, based on Monte-Carlo simulations, is that the convergence rate does not depend a priori on the dimension of the problem, and therefore should less suffer  from the curse of dimensionality encountered in deterministic procedures.

\section{Randomization approach for fully nonlinear HJB  equation}

 \setcounter{equation}{0} \setcounter{Assumption}{0}
\setcounter{Theorem}{0} \setcounter{Proposition}{0}
\setcounter{Corollary}{0} \setcounter{Lemma}{0}
\setcounter{Definition}{0} \setcounter{Remark}{0}

\subsection{Motivating example}

Let us consider the following controlled diffusion example arising from uncertain volatility model in finance: 
\beqs
X_s^{t,x,\alpha} &=& x +  \int_t^s \alpha_u dW_u, \;\;\; 0 \leq t \leq s \leq T, \;\; x \in \R, 
\enqs
where $\alpha$ is an  adapted process valued in $A$ $=$ $[\underline a,\bar a]$, $0 <  \underline a\leq \bar a < \infty$, denoted $\alpha$ $\in$ $\Ac$, interpreted as the uncertain volatility of the stock price $X$.  We define the value function of the stochastic control problem:
\beqs
v(t,x) &:=& \sup_{\alpha\in\Ac} \E\big[ h(X_T^{t,x,\alpha}) \big],
\enqs
which  is interpreted as the super-replication cost of an option payoff $h$ under uncertain volatility.  The dynamic programming equation (also called Hamilton-Jacobi-Bellman, HJB in short) 
for this stochastic control problem (see e.g. \cite{pha09}) 
is a fully nonlinear PDE in the form: 
\beq \label{Gheat}
\Dt{v} + G(D_x^2 v) &=& 0,
\enq
with terminal condition $v(T,.)$ $=$ $h$, and where 
\beqs
G(M) & := & \frac{1}{2} \sup_{a \in A}  [a^2M] \; = \; \bar a^2 M^+ - \underline{a}^2 M^-,  \;\;\; M \in \R. 
\enqs
The equation \reff{Gheat} can be viewed as a $G$-heat equation (reducing to the classical heat equation when $A$ is a singleton), and based on this observation, Peng \cite{pen06} has developed a theory of $G$-stochastic calculus with $G$-Brownian motion, extending the classical It\^o calculus with Brownian motion, and leading to the concept of nonlinear expectation. 
Denoting by $B_t^\alpha$ $=$ $\int_0^t \alpha_s dW_s$, and $\P^\alpha$ the law of $B^\alpha$ under $\P$, we notice that  $(\P^\alpha)_\alpha$ is a family of non dominated probability measures, and this contrasts with the framework of controlled drift problem in \reff{defWalpha}, which  gave rise to equivalent probability measures by Girsanov's theorem. Recalling that  Brownian motion and It\^o calculus are the basic tools for defining BSDE, Soner, Touzi and Zhang \cite{sontouzha11} have developed the theory of 2BSDE in connection with $G$-Brownian motion  by using notions  from quasi-sure analysis in a singular measures framework.  However, the main concerns with the theory of $G$-expectation and 2BSDE is that (i) it does not cover the general case of HJB equation where control appears both on the drift and diffusion, (ii) it requires uniform ellipticity condition on the diffusion coefficient, (iii) it does not lead clearly to an implementable numerical scheme since one cannot simulate a $G$-Brownian motion.  In  the rest of this paper, we shall present an alternative approach for overcoming these issues.

\subsection{BSDE with nonpositive jumps}

Let us consider the fully nonlinear PDE of HJB type:
\begin{equation} \label{HJB}
\left\{
\begin{array}{rcl}
\Dt{v} + \sup_{a\in A} \Big[ b(x,a).D_x v + \frac{1}{2}{\rm tr}(\sigma\sigma\trans(x,a) D_x^2 v)  \;\;\; & &  \\
 \;\;\;\;\;\;\;  + \;  f(x,a,v,\sigma\trans(x,a)D_x v) \Big] &=& 0, \;\;\;  \mbox{ on } \; [0,T)\times\R^d  \\
v(T,.) &=& h, \;\;\; \mbox{ on }  \;  \R^d,   
%\enq
\end{array}
\right.
\end{equation}
where $A$ is a compact metric space,  $b$ $=$ $b(x,a)$ is an $\R^d$-valued Lipschitz continuous function,  $\sigma$ $=$ $\sigma(x,a)$ is an $\R^{d\times d}$-valued  (possibly degenerate)  Lipschitz continuous function, $f$ $=$ $f(x,a,y,z)$, $h$ $=$ $h(x)$ are Lipschitz continuous functions.  Under these conditions,  there exists a unique viscosity solution with linear growth condition to \reff{HJB}, 
see \cite{ish89}. An important particular case is when $f$ $=$ $f(x,a)$ does not depend on $y,z$, and then the  PDE \reff{HJB} corresponds to the dynamic programming equation for the 
stochastic control problem: 
\beq \label{consto}
v(t,x) &=& \sup_{\alpha\in\Ac} \E \Big[ h(X_T^{t,x,\alpha}) +  \int_t^T f(X_s^{t,x,\alpha},\alpha_s) ds   \Big],
\enq
with the controlled diffusion  process in $\R^d$:
\beq \label{controlX}
X_s^{t,x,\alpha} &=& x + \int_t^s b(X_u^{t,x,\alpha},\alpha_u) du
+ \int_t^s  \sigma(X_u^{t,x,\alpha},\alpha_u) dW_u, \;\;\; t \leq s \leq T, 
\enq
where $W$ is a $d$-dimensional Brownian motion on $(\Omega,\Fc,\F,\P)$,  and $\alpha$ $\in$ $\Ac$ is the control, i.e. an $\F$-adapted process valued in $A$.  HJB type equations \reff{HJB} include the $G$-heat equation \reff{Gheat} and arise in many applications in finance, like portfolio optimization, option pricing and risk measures under model uncertainty, etc.  We refer to \cite{fleson06} or \cite{pha09} for an expository treatment of the theory of stochastic control and its applications.

The main issue for a  Feynman-Kac type formula  of the fully nonlinear PDE \reff{HJB}  comes from the fact that the controlled forward process $X^\alpha$ in \reff{controlX} cannot be simulated  for all values of the  control $\alpha$, and one cannot  remove the control process  as in the controlled drift case by Girsanov's theorem. We present here a control randomization approach, whose basic idea is to replace the control process by an (uncontrolled) auxiliary state variable process running over the control set $A$, hence simulatable, and under which one can apply Girsanov's theorem in order to recover all possible values of the original control process. As we shall see, this method allows us to provide a BSDE representation of general HJB equation \reff{HJB} in terms of a simulatable forward process formulated under a single probability measure, hence a non linear Feynman-Kac formula. An important feature of our approach is that it does not require any ellipticity condition on the diffusion coefficient. Moreover, by using a randomization with jumps, we are able to derive a practical probabilistic numerical scheme, which can take advantage of Monte-Carlo methods for dealing with high dimensional problems, both in state and control space.

\vspace{2mm}

Let us introduce a Poisson random measure $\mu(dt,da)$ on $\R_+\times A$ (hence independent of $W$), with jump times $(T_k)$, and marks $(\zeta_k)$, and intensity measure 
$\lambda(da)dt$ where $\lambda$ is a finite  measure  supporting the whole set $A$.  We denote by $\tilde\mu(dt,da)$ $=$ $\mu(dt,da)-\lambda(da)dt$ the compensated martingale measure of 
$\mu$, and  associate to $\mu$ the pure-jump process $I$ defined by:
\beqs
I_t &=& \zeta_i, \;\;\; T_k \leq t < T_{k+1}, \;\;\; k \in \N, 
\enqs
which is also written in differential form as:
\beqs
dI_t &=& \int_A (a - I_{t^-}) \mu(dt,da), \;\;\; t \geq 0.
\enqs 
We then consider the regime-switching process of dynamics:
\beqs
dX_t &=& b(X_t,I_t) dt + \sigma(X_t,I_t) dW_t, \;\;\; t \geq 0. 
\enqs 
In other words, we have replaced in the dynamics \reff{controlX}, the control $\alpha$ by the exogenous pure jump process $I$.  Notice that the pair $(X,I)$ is a Markov process valued in 
$\R^d\times A$ on the probability space $(\Omega,\Fc,\P)$ equipped  with the Brownian-Poisson filtration $\G$ $=$ $\F^{W,\mu}$ $=$ 
$(\Gc_t)_{0\leq t\leq T}$.  We next consider the BSDE with jumps, consisting in the search for a triple $(Y,Z,U)$  satisfying:
\beq
Y_t &=& h(X_T) + \int_t^T f(X_s,I_s,Y_s,Z_s) ds  \nonumber \\
& & \;\;\; - \int_t^T Z_s dW_s - \int_t^T \int_A U_s(a) \tilde\mu(ds,da), \; \;\; 0 \leq t \leq T.   \label{BSDEjumps}
\enq
Here, the pair  $(Y,Z)$ lie  in $\S_\G^2\times\L_\G^{2}(W)$,  and with respect to the Brownian framework, there is in addition the component $U$ lying in 
$\L_\G^2(\tilde\mu)$, the set of  $\Pc_\G\otimes\Bc(A)$-measurable maps  $(U_t(a))_{0\leq t\leq T}$ such that $\E[\int_0^T\int_A |U_t(a|^2 \lambda(da) dt]$ $<$ $\infty$. 
Existence and uniqueness of a triple solution $(Y,Z,U)$ $\in$ $\S_\G^2\times\L_\G^2(W)\times\L_\G^2(\tilde\mu)$ to the BSDE with jumps \reff{BSDEjumps} is proved in \cite{tanli94}, extending the result of \cite{parpen90}.   Moreover, by the Markov property of the forward regime-switching process $(X,I)$, the component solution $Y$ is written in the form:
\beqs
Y_t &=& v(t,X_t,I_t),
\enqs
for some deterministic function $v$ on $[0,T]\times\R^d\times A$,   which satisfies (in the viscosity sense) the semi-linear integro-partial differential equation (IPDE): 
\beq \label{semiIPDE}
\Dt{v} + \bar \Lc v + \Mc v + f(x,a,v,\sigma\trans(x,a)D_x v) &=& 0, \;\;\; (t,x,a) \in [0,T)\times\R^d\times A,
\enq
where 
\beqs
\bar \Lc v(t,x,a) &=& b(x,a).D_x v + \frac{1}{2}{\rm tr}(\sigma\sigma\trans(x,a)D_x^2 v), \\
\Mc v(t,x,a)  &=&  \int_A \big( v(t,x,a')- v(t,x,a) \big) \lambda(da'). 
\enqs
In other words, the BSDE with jumps  \reff{BSDEjumps} provides a Feynman-Kac formula for the semi-linear IPDE \reff{semiIPDE}, as proved in \cite{barbucpar97}, thus extending the result of \cite{parpen92}.  Moreover, when $v$ is smooth, we have by It\^o's formula, the connection:  
\beqs
Y_t  \; = \;  v(t,X_t,I_t),  & & Z_t \; = \;  \sigma\trans(X_t,I_{t^-})D_x v(t,X_t,I_{t^-}), \\
\mbox{ and } & & U_t(a) \; = \;  v(t,X_t,a) - v(t,X_t,I_{t^-}). 
\enqs

The issue is now to go from the semi-linear IPDE  \reff{semiIPDE} to the fully nonlinear PDE \reff{HJB}, and the idea is to constrain the jump component $U$ of the BSDE with jumps \reff{BSDEjumps} to be nonpositive. Let us formally derive the arguments of this approach. By constraining the jump component, this means in terms of the function $v$ (when it is smooth) that:
\beqs
U_t(a) \; = \; v(t,X_t,a) - v(t,X_t,I_{t^-}) & \leq & 0, \;\;\; \mbox{ for all } \; (t,a) \in [0,T]\times A,  
\enqs
which would imply that  $v$ $=$ $v(t,x)$ does not depend actually on $a$ $\in$ $A$.  Thus, the integral term $\Mc v$ in \reff{semiIPDE} is removed, and the variable $a$ becomes now a parameter  in the PDE satisfied by $v(t,x)$ on $[0,T)\times\R^d$:   
\beqs
\Dt{v} + b(x,a).D_x v + \frac{1}{2}{\rm tr}(\sigma\sigma\trans(x,a)D_x^2 v) +  f(x,a,v,\sigma\trans(x,a)D_x v) &=& 0, \;\;\; \mbox{ on } [0,T)\times\R^d, 
\enqs  
which should hold for any parameter value $a$ $\in$ $A$.   By taking supremum over $a$ $\in$ $A$, we formally expect to retrieve the fully nonlinear HJB equation \reff{HJB}. 
  
The rigorous derivation of the above  argument is formulated in \cite{khapha12} by means of the class of BSDE with nonpositive jumps:  this consists in a triple $(Y,Z,U)$  $\in$ 
$\S_\G^2\times\L_\G^2(W)\times\L_\G^2(\tilde\mu)$  supersolution to:
\beq
Y_t & \geq & h(X_T) + \int_t^T f(X_s,I_s,Y_s,Z_s) ds  \nonumber \\
& & \;\;\; - \int_t^T Z_s dW_s - \int_t^T \int_A U_s(a) \tilde\mu(ds,da), \; \;\; 0 \leq t \leq T.   \label{BSDEsuperjumps}
\enq
such that
\beq \label{consjump}
U_t(a) & \leq & 0\;, \;\;\;\;\; d\P\otimes dt\otimes\lambda(da) \;\;  a.e. \mbox{ on } \Omega\times[0,T]\times A. 
\enq  
By supersolution in \reff{BSDEsuperjumps}, we mean that  inequality $\geq$ holds instead of $=$ as in \reff{BSDEjumps}, and this relaxation is done   in order to get flexibility for satisfying  
the non positivity constraint in \reff{consjump}.  We are then looking for a minimal supersolution $(Y,Z,U)$ to \reff{BSDEsuperjumps}-\reff{consjump} in the sense that for any other triple 
$(\tilde Y,\tilde Z,\tilde U)$ $\in$  $\S_\G^2\times\L_\G^2(W)\times\L_\G^2(\tilde\mu)$ satisfying \reff{BSDEsuperjumps}-\reff{consjump},  we have:
\beqs
Y_t & \leq & \tilde Y_t, \;\;\; 0 \leq t \leq T, \; a.s.
\enqs
Existence and uniqueness of a minimal  supersolution $(Y,Z,U)$ to  \reff{BSDEsuperjumps}-\reff{consjump} is shown  in \cite{khapha12} by penalization methods. Moreover,  it is proved that the solution $Y$ is in the form
\beq \label{repYv}
Y_t &=& v(t,X_t),
\enq
for some deterministic function $v$ on $[0,T]\times\R^d$ (hence not depending on the state variable $I_t$, and this is the key property), and $v$ is the unique viscosity solution to the nonlinear PDE \reff{HJB}. Therefore, we have a nonlinear Feynman-Kac formula  for the solution to the fully nonlinear PDE \reff{HJB} in terms of BSDE with nonpositive jumps   \reff{BSDEsuperjumps}-\reff{consjump}.   The last paragraph shows how this probabilistic representation 
provides  a numerical scheme for solving the PDE \reff{HJB}.

\begin{Remark}
{\rm It is also proved in \cite{khapha12} that, in the case where $f$ $=$ $f(x,a)$ does not depend on $y,z$,  the minimal solution $Y$ to the BSDE with nonpositive jumps \reff{BSDEsuperjumps}-\reff{consjump} admits a dual representation in the form: 
\beqs
Y_t &=& \esssup_{\nu \in \Dc} \E^{\P^\nu} \Big[ h(X_T) + \int_t^T  f(X_s,I_s) ds \big| \Gc_t \Big], \;\;\; 0 \leq t \leq T, 
\enqs
where $\Dc$ is the set of $\Pc_\G\otimes\Bc(A)$-measurable maps $\nu$ $=$ $(\nu_t(a))_{0\leq t\leq T}$, valued in $[1,\infty)$ and bounded, and $\P^\nu$ is the probability measure  equivalent to $\P$ on $(\Omega,\Gc_T)$ whose effect by Girsanov's theorem is to change the compensator 
$\lambda(da)dt$ of $\mu$ under $\P$ to $\nu_t(a)\lambda(da)dt$ under $\P^\nu$, and to leave unchanged the Brownian motion $W$. 
Together with \reff{repYv}, this shows that the value function of a stochastic control problem \reff{consto} admits an alternative formulation in terms of intensity control. 
\ep
}
\end{Remark}

\subsection{Numerical scheme}

The discrete-time approximation  for the minimal supersolution to the BSDE with nonpositive jumps \reff{BSDEsuperjumps}-\reff{consjump} is constructed as follows. First, we observe that the pure jump process $I$ from the Poisson random measure $\mu(dt,da)$ with jump times/marks  $(T_k,\zeta_k)_k$ and intensity measure $\lambda(da)dt$ is perfectly simulated.  Indeed, the inter arrival times $S_k$ $=$ $T_{k+1}-T_k$ are i.i.d and follow an exponential law of parameter $\lambda$ $:=$ $\int_A \lambda(da)$, while the marks $\zeta_k$ are i.i.d. with distribution 
$\bar\lambda(da)$ $=$ $\lambda(da)/\lambda$, assumed to be simulatable. We thus simulate $I$  by:
\beqs
I_t &=& I_0 1_{[0,T_1)} + \sum_{k \geq 1} \zeta_k 1_{[t_k,T_{k+1})}(t), \;\;\; t \geq 0. 
\enqs
We are next given a time grid $\pi$ $:=$ $\{t_0=0 < t_1 < \ldots < t_n=T\}$ of $[0,T]$, with modulus $|\pi|$ 
$:=$ $\max_{i=1,\ldots,n} \Delta t_i$, $\Delta t_i$ $:=$ $t_{i+1}-t_i$,  and approximate the forward regime-switching process $X$ by its Euler scheme $X^\pi$ defined as
\beqs
 X^{\pi}_{t_{i+1}} &:=&  X^{\pi}_{t_i} + b(X^{\pi}_{t_i},I_{t_i}) \Delta t_i + \sigma(X^{\pi}_{t_i},I_{t_i}) \Delta W_{t_i}, \;\;\; i < n, \;\;\; X_{0}^\pi \; = \; X_0, 
 \enqs
where $\Delta W_{t_i}$ $:=$ $W_{t_{i+1}}-W_{t_i}$.  
We then  propose a discrete time approximation explicit scheme in the form: 
\begin{equation}\label{schemeBSDE}
\left\{ \begin{array}{rcl}
Y_T^\pi \; = \; \Yc_T^\pi &=& g(X_T^\pi) \\
\Zc_{t_i}^\pi &=& \E \Big[ Y_{t_{i+1}}^\pi \frac{\Delta W_{t_i}}{\Delta t_{i}} \big| \Gc_{t_i} \Big] \\
\Yc_{t_i}^\pi &=& \E\Big[ Y_{t_{i+1}}^\pi  \big| \Gc_{t_i} \Big] +   f(X_{t_i}^\pi,I_{t_i},\Yc_{t_i}^\pi,\Zc_{t_i}^\pi)  \Delta t_i \\
Y_{t_i}^\pi &=& \esssup_{a \in A} \E \Big[ \Yc_{t_i}^\pi \big| \Gc_{t_i},  I_{t_i} = a \Big], \;\;\; i = 0, \ldots, n-1.  
\end{array}
\right.
\end{equation}
The interpretation of this scheme is the following. The  first three lines in \reff{schemeBSDE} correspond to the scheme $(\Yc^\pi,\Zc^\pi)$ for a discretization of a BSDE with jumps,  as in \cite{BE08}, and exten\-ding the scheme described in paragraph 2.3  (we omit here the computation of the jump component). The last 
line in \reff{schemeBSDE} for computing the approximation $Y^\pi$ of the minimal supersolution $Y$  corresponds precisely to the minimality condition for the nonpositive jump constraint and should be understood as follows. By the Markov property of the forward process $(X,I)$, the solution $(\Yc,\Zc,\Uc)$  to the BSDE with jumps (without constraint)  is in the form $\Yc_t$ $=$ $\vartheta(t,X_t,I_t)$ for some deterministic function $\vartheta$.  Assuming that $\vartheta$ is a continuous function, 
the jump component of the BSDE, which is induced by a jump of the forward component $I$,   is equal to  $\Uc_t(a)$ $=$ $\vartheta(t,X_t,a)-\vartheta(t,X_t,I_{t^-})$.  Therefore, the nonpositive jump constraint  means that: 
$\vartheta(t,X_t,I_{t^-})$ $\geq$ $\esssup_{a\in A} \vartheta(t,X_t,a)$. The  minimality condition is thus written as: 
\beqs
Y_t \; = \;  v(t,X_t) \; = \;  \esssup_{a \in A} \vartheta(t,X_t,a) &=& \esssup_{a \in A} \E [ \Yc_t | X_t, I_{t} = a ],
\enqs
whose discrete time version is the last line in scheme \reff{schemeBSDE}.  Notice that the scheme \reff{schemeBSDE} is a dynamic programming type algorithm. The novel feature is that 
conditional expectation is taken with respect to the uncontrolled randomized extended state process $(X,I)$, and supremum with respect to the auxiliary state variable $I$.

\vspace{2mm}

The discrete-time approximation error is measured by: 
\beqs
{\rm Err}_\pm^\pi(Y) &:= &   \max_{i \leq n} 
\Big(  \E \Big[ \big(Y_{t_i} - \bar Y_{t_i}^\pi\big)_\pm^2\Big] \Big)^{1\over 2}    
%+  \sup_{t \in (t_k,t_{k+1}]} \E \Big[ \big(Y_t - \bar Y_{t_{k+1}}^\pi\big)_\pm^2 \Big] \Big)^{1\over 2} 
%+  \sup_{t \in [t_k,t_{k+1})} \E \Big[ \big(Y_t - \bar\Yc_{t_{k}}^\pi\big)_+^2 \Big] \Big)^{1\over 2} 
%\\
%{\rm Err}_-^\pi(Y) &:= &   \max_{k\leq n-1} 
%\Big(  \E \Big[ \big(Y_{t_k} - \bar Y_{t_k}^\pi\big)_-^2\Big]   +  \sup_{t \in (t_k,t_{k+1}]} \E \Big[ \big(Y_t - \bar Y_{t_{k+1}}^\pi\big)_-^2 \Big] \Big)^{1\over 2}. 
%+  \sup_{t \in [t_k,t_{k+1})} \E \Big[ \big(Y_t - \bar\Yc_{t_{k}}^\pi\big)_-^2 \Big] \Big)^{1\over 2}.
\enqs
It is proved in \cite{khalanpha13a} that 
\beqs
{\rm Err}_-^\pi(Y) & \leq & C |\pi|^{1\over 2},  
\enqs 
and under additional conditions  on $b,\sigma,f$ and $h$, namely: $b$, $\sigma$ bounded, $f$ $=$ $f(x,a,y)$ does not depend on $z$, and is convex in $y$, and $f(.,.,0)$, $g$ are bounded, 
we have
\beqs
{\rm Err}_+^\pi(Y)  & \leq &   \left\{ \begin{array}{cl}
		C |\pi|^{1\over 6} & \mbox{ when } f = f(x,a) \\
		C |\pi|^{1\over 10} & \mbox{ otherwise. } 
		\end{array}
		\right. 
\enqs
In particular, 
\beqs
  - C |\pi|^{1\over 2}  \; \leq v(0,X_0) - Y_0^\pi & \leq & \left\{ \begin{array}{cl}
		C |\pi|^{1\over 6} & \mbox{ when } f = f(x,a) \\
		C |\pi|^{1\over 10} & \mbox{ otherwise. } 
		\end{array}
		\right.
\enqs
The above error bounds are non symmetric as in deterministic methods, and are  proved by using shaking coefficients method of Krylov \cite{kry00} and switching system 
approximation of Barles and Jacobsen \cite{barjac07}.

\vspace{2mm}

The last step towards an implementable scheme consists in the approximation of the conditional expectations in  \reff{schemeBSDE}.  Here, due to the supremum operation, there is a strong advantage of using least-square regression methods. Let us briefly recall the basic principle of this method. From the definition-property of conditional expectation:
\beqs
\E[H | \Gc_{t_i} ] &=&  {\rm arg}\inf_{V \in L^2(\Gc_{t_i})} \E |H - V|^2,
\enqs
we approximate it by $\hat\E_{t_i}[H]$ $:=$ $\hat\varphi_i(X_{t_i}^\pi,I_{t_i})$ with empirical regression function:
\beqs
\hat\varphi_i &:=& {\rm arg}\inf_{\varphi \in  \Phi} \frac{1}{M} \sum_{m=1}^M \big(H^m - \varphi(X_{t_i}^{\pi,m},I_{t_i}^m) \big)^2
\enqs
where 
\begin{itemize}
\item $(X_{t_i}^{\pi,m},I_{t_i}^m)_m$ and $(H^m)_m$  are i.i.d. realizations of $(X_{t_i}^\pi,I_{t_i})$ and $H$
\item  $\Phi$ $=$  Span$\{\phi^\ell:  1 \leq \ell \leq L_\Phi\}$, $\phi^\ell$ basis functions on $\R^d\times A$.  
\end{itemize}
Then, $\esssup_{a\in A} \E[H|\Gc_{t_i},I_{t_i} =a]$  is approximated by:  
\beqs
\esssup_{a\in A}  \hat\E_{t_i,a}[ H] &:=& \esssup_{a\in A} \hat\varphi_i(X_{t_i}^\pi,a) \; = \; \hat\varphi_i(X_{t_i}^\pi, \hat a_i(X_{t_i}^\pi))
\enqs
where $\hat a_i$ is determined by:
\beqs
\hat a_i(x) &:=&  {\rm arg}\max_{a\in A} \hat\varphi_i(x,a),
\enqs
hence  in nonparametric form  given  the choice of the basis functions in $\Phi$. The advantage of this regression-projection method is that we don't need to run over the set $A$ 
in the optimization over $a$, e.g. by Newton method, as in finite-difference methods, and this quite interesting especially in high dimension for the control space $A$.  
Moreover, we  get an approximate optimal control in feedback form, i.e. a deterministic function $\hat a_i(x)$ of the state value $x$ at any date $t_i$.   
Error analysis and numerical illustrations of this algorithm are studied and performed in \cite{khalanpha13b}.

\vspace{7mm}

%\newpage 


\small 


\begin{thebibliography}{}

\bibitem{balpag03} Bally V.  and G. Pag\`es (2002):  ``A quantization algorithm for solving discrete time multidimensional optimal stopping problems", 
{\it Bernoulli},  9, 100-1049


\bibitem{barbucpar97} Barles G., R. Buckdahn, and E. Pardoux  (1997): ``Backward stochastic differential equations and integral-partial differential equations'',
{\it Stochastics and Stochastics Reports},  60, 57-83.


\bibitem{barjac07} Barles G. and E.R. Jacobsen (2007): ``Error bounds  for monotone  approximation schemes for parabolic Hamilton-Jacobi-Bellman equations", {\it Math. Computation}, 76, 1861-1893. 
 
\bibitem{bis73} Bismut J.M.  (1973): Analyse convexe et probabilit\'es, th\`ese, Facult\'e des sciences de Paris.


\bibitem{boueketou04} Bouchard B.,  Ekeland I.  and N. Touzi (2004):  ``On the Malliavin approach to Monte Carlo approximation of conditional expectations", 
{\it Finance and Stochastics},  8, 45-71.

\bibitem{BE08} Bouchard B. and R. Elie (2008):  ``Discrete time approximation of decoupled FBSDE with jumps", \textit{Stochastic Processes and Applications}, 118 (1), 53-75.


\bibitem{boutou04} Bouchard B.  and N. Touzi (2004): ``Discrete time approximation and Monte-Carlo simulation of BSDEs", {\it Stochastic processes and their applications}, 
111, 175-206.

\bibitem{elkmaz97} El Karoui, N., Mazliak L.(1997):  Backward stochastic differential equations, {\it Pitman research notes},  364.

 
\bibitem{elkpenque97}  El Karoui, N., Peng S., and M.C. Quenez (1997):  ``Backward stochastic differential equations in finance"
{\it Mathematical Finance},  7, 1-71.
 
%\bibitem{fuhpha13} Fuhrman M. and H. Pham (2013): ``Dual and backward SDE representation for optimal control of non-Markovian SDEs", to appear in {\it Annals of Applied  Probability}. 


\bibitem{fleson06} Fleming W. H. and H. M. Soner (2006):  Controlled Markov Processes and Viscosity Solutions, 2nd edition, 
Springer-Verlag.


\bibitem{gobetal06}  Gobet E.,   Lemor J.P.  and X. Warin (2006):  ``Rate of convergence of empirical regression method for
solving generalized BSDE", {\it Bernoulli},  12,  889-916. 


\bibitem{ish89} Ishii H. (1989): ``On uniqueness and existence of viscosity solutions of fully nonlinear second-order elliptic PDE's'', {\it Communications on Pure and Applied Mathematics}, 42, 15-45.

 
\bibitem{khalanpha13a} Kharroubi I., Langren\'e N. and H. Pham (2013):   ``Discrete time approximation of fully nonlinear HJB equations via BSDEs with nonpositive jumps", 
to appear in {\it Annals of Applied  Proba\-bility}. 

\bibitem{khalanpha13b} Kharroubi I., Langren\'e N. and H. Pham (2014):   ``A numerical algorithm for fully nonlinear HJB equations: an approch by control randomization",  {\it Monte-Carlo methods and applications}, 20(2), 145-165. 
 
 
%\bibitem{khaetal10} Kharroubi I., J. Ma, H. Pham, and J. Zhang (2010):  ``Backward SDEs with constrained jumps and
%Quasi-variational inequalities", {\it Annals of Probability}, 38, 794-840. 
 
 
\bibitem{khapha12}  Kharroubi I. and H. Pham (2012): ``Feynman-Kac representation for Hamilton-Jacobi-Bellman IPDE", to appear in {\it Annals of Probability}.
 
 
\bibitem{klopla92}  Kloeden P. and E. Platen (1992): Numerical Solution of Stochastic Differential Equations, Springer, Series SMAP. 

\bibitem{kry00} Krylov N. V. (2000): ``On the rate of convergence of finite difference approximations for Bellman's equations with variable coefficients",  {\it Probability Theory and Related Fields}, 117, 1-16.

\bibitem{lonsch01} Longstaff F.  and E.  Schwartz (2001):  ``Valuing american options by simulation : a simple least-square approach", 
{\it Review Of Financial Studies},  14, 113-147.

\bibitem{pagphapri04} Pag\`es G., Pham H. and J. Printems (2004): Optimal quantization and applications to numerical problems in finance, in 
{\it Handbook of Computational and Numerical Methods in Finance}, , ed. S.T. Rachev, Birkhauser, Boston. 


 
\bibitem{parpen90}  Pardoux E. and S. Peng (1990): ``Adapted solution of a backward stochastic differential equation",
{\it Systems Control Lett.},   14(1), 55-61.
 
\bibitem{parpen92}  Pardoux E. and S. Peng (1992):  ``Backward stochastic differential equations and quasilinear
parabolic partial differential equations", {\it Stochastic partial differential equations and their applications}, 
(B. L. Rozovskii and R. B. Sowers, eds.), Lect. Notes in Control and Inform. Sci., vol. 176, Springer, Berlin, pp. 200-217.
 
 
\bibitem{pen03}  Peng S. (2003): Nonlinear expectations and risk measures, in {\it Proceedings of the CIME-EMS summer school Bressanone}. 
 
\bibitem{pen06} Peng S. (2006): ``G-expectation, G-Brownian motion and related stochastic calculus of Ito type",
Proceedings of 2005, Abel symposium, Springer. 
 
\bibitem{pha09} Pham H. (2009): Continuous time stochastic control and optimization with financial applications, Springer, Series SMAP. 


\bibitem{sontouzha11} Soner M., Touzi N.  and J. Zhang (2011): ``The wellposedness of second order backward SDEs",
{\it Probability Theory and Related Fields}, 153, 149-190. 


\bibitem{tanli94} Tang S. and X. Li (1994): ``Necessary conditions for optimal control of stochastic systems with jumps'', {\it SIAM J. Control and Optimization}, 32, 1447-1475.



\bibitem{zha04} Zhang J. (2004): ``A numerical scheme for BSDEs", {\it Annals of Applied Probability}, 14, 459-488.
 
 
 

\end{thebibliography}
\end{document}